\documentclass{article}
\usepackage{amsfonts, amsmath, amssymb}

\begin{document}

\begin{center}
{\Large Remarks regarding some special matrices}%
\[
\]

Cristina FLAUT and Andreea BAIAS

\[
\]
\end{center}

\textbf{Abstract. }In this paper, by using matix representation for
quaternions and octonions, we provide a procedure to obtain some example of $%
k-$potent matrices of order $4$ or $8$, over the real field or over the
field $\mathbb{Z}_{p}$, with $p$ a prime number.  
\[
\]

\textbf{1. Introduction}%
\[
\]

In this paper, in the study of $k-$potent elements we extend the results
obtained in [FB; 24] to generalised quaternion algebras and generalised
octonion algebras. As an application, by using the matrix representations
for quaternions and octonions, we give examples of classes of matrices which
are $k-$potents and a procedure  to find such examples.

The paper is organised as follows: in the second paragraph we\ present
matrix representations and their properties which can be used to provide
examples of $k-$potent matrices and in the part three we study $k-$potents
elements in quaternion algebras and octonion algebras over the real field.
The paper end with conclusions and an idea for a further research. 

\[
\]%
\textbf{2.} \textbf{Matrix representation}%
\[
\]

In this paper, the field $K$ \ is considered with characteristic different
from two. In the following, we will consider the quaternion algebra over an
arbitrary field $K$. 

For \ two elements $a,b\in K$, we define a generalized quaternion algebra,
denoted by $\mathbb{H(}\alpha ,\beta \mathbb{)}=\left( \frac{a,b}{K}\right) $%
, with basis $\{1,f_{1},f_{2},f_{3}\}$ and multiplication given in the
following table:%
\[
\begin{tabular}{l|llll}
$\cdot $ & $1$ & $f_{1}$ & $f_{2}$ & $f_{3}$ \\ \hline
$1$ & $1$ & $f_{1}$ & $f_{2}$ & $f_{3}$ \\ 
$f_{1}$ & $f_{1}$ & $a$ & $f_{3}$ & $af_{2}$ \\ 
$f_{2}$ & $f_{2}$ & $-f_{3}$ & $b$ & $-bf_{1}$ \\ 
$f_{3}$ & $f_{3}$ & $-af_{2}$ & $bf_{1}$ & $-ab$%
\end{tabular}%
\]

If $q\in \mathbb{H}(a,b),$ $q=q_{0}+q_{1}f_{1}+q_{2}f_{2}+q_{3}f_{3},$ then 
\[
\overline{q}=q_{0}-q_{1}f_{1}-q_{2}f_{2}-q_{3}f_{3} 
\]%
is called the \textit{conjugate} of the element $q.$ For $q\in \mathbb{H}%
(a,b),$ we consider the following elements:

\[
\mathbf{t}\left( q\right) =q+\overline{q}\in K
\]%
and

\[
\,\mathbf{n}\left( q\right) =q\overline{q}%
=q_{0}^{2}-aq_{1}^{2}-bq_{2}^{2}+abq_{3}^{2}\in K, 
\]%
called the \textit{trace}, respectively, the \textit{norm} of the element $%
q\in \mathbb{H}(a,b)$. \thinspace \thinspace It\thinspace \thinspace
\thinspace follows\thinspace \thinspace \thinspace that$\,\,$%
\[
\left( q+\overline{q}\right) q\,=q^{2}+\overline{q}q=q^{2}+\mathbf{n}\left(
q\right) \cdot 1 
\]%
and\thinspace \thinspace 
\[
q^{2}-\mathbf{t}\left( q\right) q+\mathbf{n}\left( q\right) =0,\forall q\in 
\mathbb{H}(a,b),\, 
\]%
therefore the generalized quaternion algebra is a \textit{quadratic algebra}.

If, \thinspace for $x\in \mathbb{H}(a,b),$ the relation $\mathbf{n}\left(
x\right) =0$ implies $x=0$, then the algebra $\mathbb{H}(a,b)$ is a \textit{%
division} algebra. A quaternion non-division algebra is called a \textit{%
split} algebra.

Using the above notations, we remark that $\mathbb{H}\left( -1,-1\right)
=\left( \frac{-1,-1}{\mathbb{R}}\right) $ is a division algebra.

A generalized octonion algebra over an arbitrary field $K,$ with \textit{char%
}$K\neq 2,$ is an algebra of dimension $8,$ denoted $\mathbb{O}(a,b,c),$
with basis $\{1,f_{1},...,f_{7}\}\medskip \,\ $and multiplication given in
the following table:\medskip \medskip 

\begin{center}
{\footnotesize $%
\begin{tabular}{c||c|c|c|c|c|c|c|c|}
$\cdot $ & $1$ & $\,\,\,f_{1}$ & $\,\,\,\,\,f_{2}$ & $\,\,\,\,f_{3}$ & $%
\,\,\,\,f_{4}$ & $\,\,\,\,\,\,f_{5}$ & $\,\,\,\,\,\,f_{6}$ & $%
\,\,\,\,\,\,\,f_{7}$ \\ \hline\hline
$\,1$ & $1$ & $\,\,\,f_{1}$ & $\,\,\,\,f_{2}$ & $\,\,\,\,f_{3}$ & $%
\,\,\,\,f_{4}$ & $\,\,\,\,\,\,f_{5}$ & $\,\,\,\,\,f_{6}$ & $%
\,\,\,\,\,\,\,f_{7}$ \\ \hline
$\,f_{1}$ & $\,\,f_{1}$ & $\alpha $ & $\,\,\,\,f_{3}$ & $af_{2}$ & $%
\,\,\,\,f_{5}$ & $af_{4}$ & $-\,\,f_{7}$ & $\,\,-af_{6}$ \\ \hline
$\,f_{2}$ & $\,f_{2}$ & $-f_{3}$ & $b$ & $\,\,-bf_{1}$ & $\,\,\,\,f_{6}$ & $%
\,\,\,\,\,f_{7}$ & $bf_{4}$ & $bf_{5}$ \\ \hline
$f_{3}$ & $f_{3}$ & -$af_{2}$ & $bf_{1}$ & $-ab$ & $\,\,\,\,f_{7}$ & $af_{6}$
& $\,\,\,-bf_{5}$ & $-abf_{4}$ \\ \hline
$f_{4}$ & $f_{4}$ & $-f_{5}$ & $-\,f_{6}$ & $-\,\,f_{7}$ & $\,c$ & $%
\,\,-\,cf_{1}$ & $\,\,-cf_{2}$ & $\,\,\,\,-\,cf_{3}$ \\ \hline
$\,f_{5}$ & $\,f_{5}$ & -$af_{4}$ & $-\,f_{7}$ & -$\,af_{6}$ & $cf_{1}$ & -$%
ac$ & $cf_{3}$ & $\,acf_{2}$ \\ \hline
$\,\,f_{6}$ & $\,\,f_{6}$ & $\,\,\,\,f_{7}$ & $\,\,-bf_{4}$ & $bf_{5}$ & $%
cf_{2}$ & $\,\,\,-cf_{3}$ & -$bc$ & $-bcf_{1}$ \\ \hline
$\,\,f_{7}$ & $\,\,f_{7}$ & $af_{6}$ & $\,-bf_{5}$ & $abf_{4}$ & $cf_{3}$ & $%
-acf_{2}$ & $bcf_{1}$ & $abc$ \\ \hline
\end{tabular}%
\ \medskip $ }\medskip 
\end{center}

The algebra $\mathbb{O}(a,b,c)$ is a non-commutative and a non-associative
algebra, but it is \textit{alternative},\thinspace \textit{\ flexible} and 
\textit{power-associative}.

If $x\in \mathbb{O}(a,b,c),$ $%
x=x_{0}+x_{1}f_{1}+x_{2}f_{2}+x_{3}f_{3}+x_{4}f_{4}+x_{5}f_{5}+x_{6}f_{6}+x_{7}f_{7},
$ then $\overline{x}%
=x_{0}-x_{1}f_{1}-x_{2}f_{2}-x_{3}f_{3}-x_{4}f_{4}-x_{5}f_{5}-x_{6}f_{6}-x_{7}f_{7}
$ is called the \textit{conjugate} of the element $x.$ For $x\in \mathbb{O}%
(a,b,c),$ we define the elements:

\[
\mathbf{t}\left( x\right) =x+\overline{x}\in K
\]%
and 
\[
\,\mathbf{n}\left( x\right) =x\overline{x}%
=x_{0}^{2}-ax_{1}^{2}-bx_{2}^{2}+abx_{3}^{2}-cx_{4}^{2}+acx_{5}^{2}+bcx_{6}^{2}-abcx_{7}^{2}\in K.
\]%
These elements are called the \textit{trace}, respectively, the \textit{norm}
of the element $x\in \mathbb{O}(a,b,c).$ \thinspace \thinspace It\thinspace
\thinspace \thinspace follows\thinspace \thinspace \thinspace that$\,\,$%
\[
\left( x+\overline{x}\right) x\,=x^{2}+\overline{x}x=x^{2}+\mathbf{n}\left(
x\right) \cdot 1
\]%
and\thinspace \thinspace 
\[
x^{2}-\mathbf{t}\left( x\right) x+\mathbf{n}\left( x\right) =0,\forall x\in
x\in \mathbb{O}(a,b,c),
\]%
$\,$therefore the generalized octonion algebra is a \textit{quadratic }%
algebra.

If, \thinspace for $x\in \mathbb{O}(a,b,c),$ the relation $\mathbf{n}\left(
x\right) =0$ implies $x=0$, then the algebra $\mathbb{O}(a,b,c)$ is a 
\textit{division} algebra ( see [Sc; 54] and [Sc; 66]).

If we take $a=b=c=-1$, $K=\mathbb{R}$, then we obtain $\mathbb{H}\left(
-1,-1\right) $, the quaternion division algebra, usually denoted by $\mathbb{%
H}$ and octonion division algebra $\mathbb{O}\left( -1,-1,-1\right) $,
usually denoted by $\mathbb{O}$. For example, by taking $a=-1$ and $b=1$, $K=%
\mathbb{R}$, we obtain a split quaternion algebra. In the following, we will
denote $\mathbb{H}_{K}=\left( \frac{-1,-1}{K}\right) $ and $\mathbb{O}%
_{K}=\left( \frac{-1,-1,-1}{K}\right) $.

We know that a finite-dimensional associative algebra $A$ over an arbitrary
field $K$ is algebraically isomorphic to a subalgebra of a matrix algebra
over the same field $K$. Therefore, each element $a\in A$ has a matrix
representation. That means, there is a map $f:A\rightarrow \mathcal{M}%
_{n}\left( K\right) $ such that $f\left( x\right) =M_{x}\in \mathcal{M}%
_{n}\left( K\right) $, where \textit{dim}$A=n$. For an arbitrary quaternion
algebra $\mathbb{H}\left( a,b\right) $, the map 
\[
\varphi :\mathbb{H}\left( a,b\right) \rightarrow \mathcal{M}_{4}\left(
K\right) 
\]%
\begin{equation}
\varphi \left( q\right) =\left( 
\begin{array}{cccc}
q_{0} & aq_{1} & bq_{2} & -abq_{3} \\ 
q_{1} & q_{0} & bq_{3} & -bq_{2} \\ 
q_{2} & -aq_{3} & q_{0} & aq_{1} \\ 
q_{3} & -q_{2} & q_{1} & q_{0}%
\end{array}%
\right) ,q\in \mathbb{H}\left( a,b\right) \text{,}  \tag{1}
\end{equation}%
$\allowbreak $is called \textit{the left representation} and the map

\[
\rho :\mathbb{H}\left( a,b\right) \rightarrow \mathcal{M}_{4}\left( K\right) 
\]%
\begin{equation}
\rho \left( q\right) =\left( 
\begin{array}{cccc}
q_{0} & aq_{1} & bq_{2} & -abq_{3} \\ 
q_{1} & q_{0} & -bq_{3} & bq_{2} \\ 
q_{2} & aq_{3} & q_{0} & -aq_{1} \\ 
q_{3} & q_{2} & -q_{1} & q_{0}%
\end{array}%
\right) ,q\in \mathbb{H}\left( a,b\right) \text{,}  \tag{2}
\end{equation}%
where $q=q_{0}+q_{1}f_{1}+q_{2}f_{2}+q_{3}f_{3}$, is called \textit{the
right representation} ( see [Ti; 00]).

In the same paper [Ti; 00], were defined, for real octonions, two
representations maps, left and right representations, by using the maps $%
\varphi $ and $\rho $, defined on real quaternions. These maps can be
defined for all octonion algebras $\mathbb{O}\left( a,b,c\right) $ over an
arbitrary field \thinspace $K$, namely%
\[
\Phi :\mathbb{O}\left( a,b,c\right) \rightarrow \mathcal{M}_{8}\left(
K\right) , 
\]%
\[
\Phi \left( x\right) =\left( 
\begin{array}{cc}
\varphi \left( x^{\prime }\right) & -\rho \left( x^{\prime \prime }\right)
E_{4} \\ 
\varphi \left( x^{\prime \prime }\right) E_{4} & \rho (x^{\prime })%
\end{array}%
\right) , 
\]%
the left representation, where the octonion $x$ can be written under the
form $x=x^{\prime }+x^{\prime \prime }f\,,$ with $x^{\prime },x^{\prime
\prime }\in \mathbb{H}\left( a,b\right) $, by using the Cayley-Dickson
process. In the same way, we define the right representation 
\[
\Psi :\mathbb{O}\left( a,b,c\right) \rightarrow \mathcal{M}_{8}\left(
K\right) , 
\]%
\[
\Psi \left( x\right) =\left( 
\begin{array}{cc}
\rho \left( x^{\prime }\right) & -\varphi \left( \overline{x^{\prime \prime }%
}\right) \\ 
\varphi \left( x^{\prime }\right) & \rho (\overline{x^{\prime }})%
\end{array}%
\right) , 
\]%
where $\overline{x^{\prime }},$ $\overline{x^{\prime \prime }}$ are the
conjugates of the quaternions $x^{\prime }$ and $x^{\prime \prime }$ and $%
E_{4}=diag\left( 1,-1,-1,-1\right) $.\smallskip\ 

For $x\in \mathbb{O}\left( a,b,c\right) $, $%
x=x_{0}+x_{1}f_{1}+x_{2}f_{2}+x_{3}f_{3}+x_{4}f_{4}+x_{5}f_{5}+x_{6}f_{6}+x_{7}f_{7}
$ with $\{1,f_{1},....,f_{7}\}$ the base in $\mathbb{O}\left( a,b,c\right) $%
, we have 
\begin{equation}
\Phi \left( x\right) =\left( 
\begin{array}{cccccccc}
x_{0} & ax_{1} & bx_{2} & -abx_{3} & cx_{4} & -acx_{5} & -bcx_{6} & abcx_{7}
\\ 
x_{1} & x_{0} & bx_{3} & -bx_{2} & cx_{5} & -cx_{4} & bcx_{7} & -bcx_{6} \\ 
x_{2} & -ax_{3} & x_{0} & ax_{1} & cx_{6} & -acx_{7} & -cx_{4} & acx_{5} \\ 
x_{3} & -x_{2} & x_{1} & x_{0} & cx_{7} & -cx_{6} & cx_{5} & -cx_{4} \\ 
x_{4} & -ax_{5} & -bx_{6} & abx_{7} & x_{0} & ax_{1} & bx_{2} & -abx_{3} \\ 
x_{5} & -x_{4} & -bx_{7} & bx_{6} & x_{1} & x_{0} & -bx_{3} & bx_{2} \\ 
x_{6} & ax_{7} & -x_{4} & -ax_{5} & x_{2} & ax_{3} & x_{0} & -ax_{1} \\ 
x_{7} & x_{6} & -x_{5} & -x_{4} & x_{3} & x_{2} & -x_{1} & x_{0}%
\end{array}%
\right)   \tag{3}
\end{equation}%
and 
\begin{equation}
\Psi \left( x\right) =\left( 
\begin{array}{cccccccc}
x_{0} & ax_{1} & bx_{2} & -abx_{3} & cx_{4} & -acx_{5} & -bcx_{6} & abcx_{7}
\\ 
x_{1} & x_{0} & -bx_{3} & bx_{2} & -cx_{5} & cx_{4} & -bcx_{7} & bcx_{6} \\ 
x_{2} & ax_{3} & x_{0} & -ax_{1} & -cx_{6} & acx_{7} & cx_{4} & -acx_{5} \\ 
x_{3} & x_{2} & -x_{1} & x_{0} & -cx_{7} & cx_{6} & -cx_{5} & cx_{4} \\ 
x_{4} & ax_{5} & bx_{6} & -abx_{7} & x_{0} & -ax_{1} & -bx_{2} & abx_{3} \\ 
x_{5} & x_{4} & bx_{7} & -bx_{6} & -x_{1} & x_{0} & bx_{3} & -bx_{2} \\ 
x_{6} & -ax_{7} & x_{4} & ax_{5} & -x_{2} & -ax_{3} & x_{0} & ax_{1} \\ 
x_{7} & -x_{6} & x_{5} & x_{4} & -x_{3} & -x_{2} & x_{1} & x_{0}%
\end{array}%
\right)   \tag{4}
\end{equation}

\textbf{Proposition 1.} ([Ti; 00], Lemma 1.2) \textit{With the above
notations, for} $\varepsilon \in \{\varphi ,\rho \}$\textit{, we have:}

i) $\varepsilon \left( x+y\right) =\varepsilon \left( x\right) +\varepsilon
\left( y\right) ,$

ii) $\varepsilon \left( xy\right) =\varepsilon \left( x\right) \varepsilon
\left( y\right) ,$

iii) $\varepsilon \left( \lambda x\right) =\lambda \varepsilon \left(
x\right) ,\lambda \in K,\varepsilon \left( 1\right) =I_{4}$

iv) $\varepsilon \left( \overline{x}\right) =\varepsilon ^{T}\left( x\right) 
$

v) $\varepsilon \left( x^{-1}\right) =\varepsilon ^{-1}\left( x\right) $

vi) $\varepsilon \left( x\right) =\varepsilon \left( y\right) $ if and only
if $x=y$, with $x,y\in \mathbb{H}$.$\smallskip \smallskip $

We remark that the above properties are proved in the real case for division
algebra $\mathbb{H}$, but these are true in the general case, for an
arbitrary field with characteristic different from $2$. Properties i), iii),
iv) and v) are also satisfied for the maps $\Phi $ and $\Psi $ over reals.
Moreover, the following properties were proved over the real
field.\smallskip 

\textbf{Proposition 2.} ([Ti;00],Theorem 2.5, Theorem 2.9, Theorem 2,10,
Theorem 2.11)

\textit{With the above notations, for} $\varepsilon \in \{\Phi ,\Psi \}$%
\textit{, we have}

i) $\varepsilon \left( x+y\right) =\varepsilon \left( x\right) +\varepsilon
\left( y\right) ,$

ii) $\varepsilon \left( x^{2}\right) =\varepsilon \left( x\right) ^{2}$

iii) $\varepsilon \left( xyx\right) =\varepsilon \left( x\right) \varepsilon
\left( y\right) \varepsilon \left( x\right) ,$

iv) $\varepsilon \left( \lambda x\right) =\lambda \varepsilon \left(
x\right) ,\lambda \in K,\varepsilon \left( 1\right) =I_{8}$

v) $\varepsilon \left( \overline{x}\right) =\varepsilon ^{T}\left( x\right) $

vi)$~\varepsilon \left( x\right) =\varepsilon \left( y\right) $ if and only
if $x=y$, where $x,y\in \mathbb{O}$.\smallskip 

Properties ii) and iii) from the above proposition were proved by using
alternativity and \ the the Moufang identities. But, since Moufang
identities are true in any octonion algebra over a field of characteristic
not two (see \textbf{[}Sc; 66\textbf{]}), these properties are also true for 
$K=\mathbb{Z}_{p}$, $p$ a prime number, $p\neq 2$.\smallskip

\textbf{Proposition 3.} \textit{With the notations from the above
proposition, we have that} 
\[
\varepsilon \left( x^{n}\right) =\varepsilon ^{n}\left( x\right) \text{,}
\]%
\textit{for} $x\in \mathbb{O}\left( a,b,c\right) $ \textit{over an arbitrary
field} $K$ \textit{and} $n$ \textit{a positive integer}.\smallskip 

\textbf{Proof.} We use induction. From condition ii) from the above
proposition, taking $x=y$, we have \ $\varepsilon \left( x^{3}\right)
=\varepsilon ^{3}\left( x\right) $. From the same condition, for $y=x^{2}$,
we have $\varepsilon \left( x^{4}\right) =\varepsilon ^{4}\left( x\right) $
and so on.\smallskip\ 

In the paper [FB; 24], we studied some properties of $k-$potent elements
over algebras obtained by the Cayley-Dickson process.\smallskip 

\textbf{Definition 4.}

i) The element $x$ in the ring $R$ is called \textit{nilpotent }if there is
a positive integer $n$ such that $x^{n}=0$. The number $n$ is the smallest
with this property and is called the \textit{nilpotency index}. \ ii) The
element $x$ in the ring $R$ is called a $k$\textit{-potent} element, for $%
k>1 $, a positive integer, if $k$ is the smallest number such that $x^{k}=x$%
. The number $k$ is called the $k$\textit{-potency index}. For $k=2,$ we
have idempotent elements, for $k=3$, we have tripotent elements, and so
on.\smallskip \smallskip

\textbf{Remark 5. }From the above definition and Proposition 3, if $x\in 
\mathbb{H}\left( a,b\right) $ or $x\in \mathbb{O}\left( a,b,c\right) $ is a $%
k-$potent element, it results that the matrices $\varepsilon \left( x\right)
,$ for $\varepsilon \in \{\varphi ,\rho \}$ or $\varepsilon \in \{\Phi ,\Psi
\}$ are $k-$potent matrices over the field $K$.\smallskip 

\textbf{Example 6.} i) By using results obtained in [FB; 24], we consider
quaternions over the field\textbf{\ }$K=\mathbb{Z}_{5}$ and the element%
\textbf{\ }$x=2+3i+j+3k$ which is a $5-$potent element over $\mathbb{H}_{%
\mathbb{Z}_{5}}$. Indeed, $x=2+3\gamma ,\gamma =i+2j+k$, with $\gamma ^{2}=-1
$ and $\mathbf{n}_{x}=1$. Therefore, $x^{2}=2\gamma $ and $x^{4}=1$.  The
matrices 
\[
\varphi \left( x\right) =\left( 
\begin{array}{cccc}
2 & -3 & -1 & -3 \\ 
3 & 2 & -3 & 1 \\ 
1 & 3 & 2 & -3 \\ 
3 & -1 & 3 & 2%
\end{array}%
\right) \text{ and }\rho \left( x\right) =\left( 
\begin{array}{cccc}
2 & -3 & -1 & -3 \\ 
3 & 2 & 3 & -1 \\ 
1 & -3 & 2 & 3 \\ 
3 & 1 & -3 & 2%
\end{array}%
\right) 
\]%
are $5-$potent matrices, that means $\varphi ^{4}\left( z\right) =\rho
^{4}\left( z\right) =I_{4}$.

ii) With the same arguments as above, if we consider octonions over the field%
\textbf{\ }$K=\mathbb{Z}_{13}$ and the element $x\in \mathbb{O}_{\mathbb{Z}%
_{13}}$, $x=3+2f_{1}+f_{2}+f_{3}+f_{4}+f_{5}+f_{6}+f_{7}$, we obtain that $x$
is a $13-$potent element . The matrices 
\[
\Phi \left( x\right) =\left( 
\begin{array}{cccccccc}
3 & -2 & -1 & -1 & -1 & -1 & -1 & -1 \\ 
2 & 3 & -1 & 1 & -1 & 1 & 1 & -1 \\ 
1 & 1 & 3 & -2 & -1 & -1 & 1 & 1 \\ 
1 & -1 & 2 & 3 & -1 & 1 & -1 & 1 \\ 
1 & 1 & 1 & 1 & 3 & -2 & -1 & -1 \\ 
1 & -1 & 1 & -1 & 2 & 3 & 1 & -1 \\ 
1 & -1 & -1 & 1 & 1 & -1 & 3 & 2 \\ 
1 & 1 & -1 & -1 & 1 & 1 & -2 & 3%
\end{array}%
\right) 
\]%
and 
\[
\Psi \left( x\right) =\left( 
\begin{array}{cccccccc}
3 & -2 & -1 & -1 & -1 & -1 & -1 & -1 \\ 
2 & 3 & 1 & -1 & 1 & -1 & -1 & 1 \\ 
1 & -1 & 3 & 2 & 1 & 1 & -1 & -1 \\ 
1 & 1 & -2 & 3 & 1 & -1 & 1 & -1 \\ 
1 & -1 & -1 & -1 & 3 & 2 & 1 & 1 \\ 
1 & 1 & -1 & 1 & -2 & 3 & -1 & 1 \\ 
1 & 1 & 1 & -1 & -1 & 1 & 3 & -2 \\ 
1 & -1 & 1 & 1 & -1 & -1 & 2 & 3%
\end{array}%
\right) 
\]%
are also $13$-potent matrices.%
\[
\]%
\textbf{3.} $\mathbf{k-}$\textbf{potents elements in quaternion algebras and
octonion algebras over the real field}%
\[
\]

In the following, we will consider $\mathbb{H}$ and $\mathbb{O}$ the real
division quaternion algebra and the real division octonion algebra. Let $%
\mathbb{A}\in \{\mathbb{H}$, $\mathbb{O}\}$. If $x\in \mathbb{A}$. $x\neq 0$%
, is a $k-$potent element, therefore $x^{k}=x$. Since $\mathbb{A}$ is a
division algebra, we have that $x$ is an invertible element, therefore $%
n_{x}\neq 0$ and $x^{k-1}=1$. It is clear from here that a $k-$potent
element is a solution of the equation 
\begin{equation}
x^{k-1}=1\text{.}  \tag{5}
\end{equation}%
In the following, we will provide solutions for this equation. From relation 
$\left( 5\right) $, we obtain that $n_{x}^{k-1}=1$. Since $n_{x}$ is a
positive number, we obtain that $n_{x}=1$. Let $x\in \mathbb{H}$, $%
x=x_{0}+x_{1}f_{1}+x_{2}f_{2}+x_{3}f_{3}$, with $%
n_{x}=x_{0}^{2}+x_{1}^{2}+x_{2}^{2}+x_{3}^{2}=1,$ with $%
x_{0},x_{1},x_{2},x_{3}\in (-1,1)$. We denote with 
\[
\cos \alpha =x_{0},\sin \alpha =\sqrt{x_{1}^{2}+x_{2}^{2}+x_{3}^{2}},
\]%
\[
\theta =\frac{x_{1}f_{1}+x_{2}f_{2}+x_{3}f_{3}}{\sqrt{%
x_{1}^{2}+x_{2}^{2}+x_{3}^{2}}}\in \mathbb{H}\text{, with }\theta ^{2}=-1.
\]%
and the element $x$ can be write under the form%
\[
x=\cos \alpha +\theta \sin \alpha \text{. }
\]%
It is clear that $x^{2}=\left( \cos \alpha +\theta \sin \alpha \right)
^{2}=\cos ^{2}\alpha -\sin ^{2}\theta +(2\cos \alpha \sin \alpha )\theta
=\cos 2\alpha +\theta \sin 2\alpha $. Using induction, we obtain that 
\[
x^{n}=\cos n\alpha +\theta \sin n\alpha .
\]%
By using the above relation, the element $x\in \mathbb{H}$ satisfing
condition $x^{k-1}=1$ has the form 
\[
x=\cos \frac{2\pi }{k-1}+\theta \sin \frac{2\pi }{k-1}\text{.}
\]

\textbf{Example 7.} i) We consider $x=\frac{1}{2}+\frac{1}{2}f_{1}+\frac{1}{2%
}f_{2}+\frac{1}{2}f_{3}=\cos \frac{\pi }{3}+\theta \sin \frac{\pi }{3}$,
where $\theta =\frac{f_{1}+f_{2}+f_{3}}{\sqrt{3}}$. We obtain that $x^{6}=1$%
, then $x^{7}=x$ and $x$ is a $7-$potent element. Therefore, the matrices 
\[
\varphi \left( x\right) =\left( 
\begin{array}{cccc}
\frac{1}{2} & -\frac{1}{2} & -\frac{1}{2} & -\frac{1}{2} \\ 
\frac{1}{2} & \frac{1}{2} & -\frac{1}{2} & \frac{1}{2} \\ 
\frac{1}{2} & \frac{1}{2} & \frac{1}{2} & -\frac{1}{2} \\ 
\frac{1}{2} & -\frac{1}{2} & \frac{1}{2} & \frac{1}{2}%
\end{array}%
\right) ~\text{and~}\rho \left( x\right) =\left( 
\begin{array}{cccc}
\frac{1}{2} & -\frac{1}{2} & -\frac{1}{2} & -\frac{1}{2} \\ 
\frac{1}{2} & \frac{1}{2} & \frac{1}{2} & -\frac{1}{2} \\ 
\frac{1}{2} & -\frac{1}{2} & \frac{1}{2} & \frac{1}{2} \\ 
\frac{1}{2} & \frac{1}{2} & -\frac{1}{2} & \frac{1}{2}%
\end{array}%
\right) 
\]%
are $7-$potent matrices.

ii) For $x=-\frac{1}{2}+\frac{1}{2}f_{1}-\frac{1}{2}f_{2}+\frac{1}{2}%
f_{3}=\cos \frac{2\pi }{3}+\theta \sin \frac{2\pi }{3}$, where $\theta =%
\frac{f_{1}-f_{2}+f_{3}}{\sqrt{3}}$. We obtain that $x^{3}=1,$ then $x^{4}=x$
and $x$ is $4-$potent element. The matrices

\[
\varphi \left( q\right) =\left( 
\begin{array}{cccc}
-\frac{1}{2} & -\frac{1}{2} & \frac{1}{2} & -\frac{1}{2} \\ 
\frac{1}{2} & -\frac{1}{2} & -\frac{1}{2} & -\frac{1}{2} \\ 
-\frac{1}{2} & \frac{1}{2} & -\frac{1}{2} & -\frac{1}{2} \\ 
\frac{1}{2} & \frac{1}{2} & \frac{1}{2} & -\frac{1}{2}%
\end{array}%
\right) \text{ and }\rho \left( q\right) =\left( 
\begin{array}{cccc}
-\frac{1}{2} & -\frac{1}{2} & \frac{1}{2} & -\frac{1}{2} \\ 
\frac{1}{2} & -\frac{1}{2} & \frac{1}{2} & \frac{1}{2} \\ 
-\frac{1}{2} & -\frac{1}{2} & -\frac{1}{2} & \frac{1}{2} \\ 
\frac{1}{2} & -\frac{1}{2} & -\frac{1}{2} & -\frac{1}{2}%
\end{array}%
\right) \text{ } 
\]%
are $4-$potent matrices.

iii) For $x=\frac{1}{2}f_{1}-\frac{1}{2}f_{2}+\frac{\sqrt{2}}{2}f_{3}=\cos 
\frac{\pi }{2}+\theta \sin \frac{\pi }{2}$, where $\theta =\frac{f_{1}-f_{2}+%
\sqrt{2}f_{3}}{2}$. We obtain that $x^{4}=1,$ then $x^{5}=x$ and $x$ is $5-$%
potent element. The matrices 
\[
\varphi \left( q\right) =\left( 
\begin{array}{cccc}
0 & -\frac{1}{2} & \frac{1}{2} & -\frac{\sqrt{2}}{2} \\ 
\frac{1}{2} & 0 & -\frac{\sqrt{2}}{2} & -\frac{1}{2} \\ 
-\frac{1}{2} & \frac{\sqrt{2}}{2} & 0 & -\frac{1}{2} \\ 
\frac{\sqrt{2}}{2} & \frac{1}{2} & \frac{1}{2} & 0%
\end{array}%
\right) \text{ and }\rho \left( q\right) =\left( 
\begin{array}{cccc}
0 & -\frac{1}{2} & \frac{1}{2} & -\frac{\sqrt{2}}{2} \\ 
\frac{1}{2} & 0 & \frac{\sqrt{2}}{2} & \frac{1}{2} \\ 
-\frac{1}{2} & -\frac{\sqrt{2}}{2} & 0 & \frac{1}{2} \\ 
\frac{\sqrt{2}}{2} & -\frac{1}{2} & -\frac{1}{2} & 0%
\end{array}%
\right) \text{ } 
\]%
are $5-$potent matrices.

\textbf{Remark 8.} In the following, we consider the quaternion split
algebra $\mathbb{H}\left( a,b\right) $ or octonion split algebra $\mathbb{O}%
\left( a,b,c\right) $ and $K=\mathbb{R}$. If an element $x\in \mathbb{H}%
\left( a,b\right) $ or $x\in \mathbb{O}\left( a,b,c\right) $ is $k-$potent
with $n_{x}=0$ and $t_{x}\neq 0$, it results that $t_{x}^{k}=t_{x}$,
therefore $t_{x}=1$ or $t_{x}=-1$ and $x_{0}=\frac{t_{x}}{2}$. We have that $%
t_{x}^{2}=t_{x}$ or $t_{x}^{3}=t_{x}$ , therefore there are only idempotent
and tripotent elements in $\mathbb{H}\left( a,b\right) $ or $\mathbb{O}%
\left( a,b,c\right) $ over $\mathbb{R}$.\smallskip 

\textbf{Example 9.} 

i) We consider the split quaternion algebra $\mathbb{H}\left( 1,1\right) $
and the quaternion $q=\frac{1}{2}\left( 1+f_{1}+f_{2}+f_{3}\right) $, with $%
\mathbf{n}\left( q\right) =0$. We have that the matrices 
\[
\varphi \left( q\right) =\left( 
\begin{array}{cccc}
\frac{1}{2} & \frac{1}{2} & \frac{1}{2} & -\frac{1}{2} \\ 
\frac{1}{2} & \frac{1}{2} & \frac{1}{2} & -\frac{1}{2} \\ 
\frac{1}{2} & -\frac{1}{2} & \frac{1}{2} & \frac{1}{2} \\ 
\frac{1}{2} & -\frac{1}{2} & \frac{1}{2} & \frac{1}{2}%
\end{array}%
\right) \text{ and }\rho \left( q\right) =\left( 
\begin{array}{cccc}
\frac{1}{2} & \frac{1}{2} & \frac{1}{2} & -\frac{1}{2} \\ 
\frac{1}{2} & \frac{1}{2} & -\frac{1}{2} & \frac{1}{2} \\ 
\frac{1}{2} & \frac{1}{2} & \frac{1}{2} & -\frac{1}{2} \\ 
\frac{1}{2} & \frac{1}{2} & -\frac{1}{2} & \frac{1}{2}%
\end{array}%
\right) ,
\]%
are idempotent matrices. For $w=\frac{1}{2}\left(
-1+f_{1}+f_{2}+f_{3}\right) $, the matrices 

\[
\varphi \left( w\right) =\left( 
\begin{array}{cccc}
-\frac{1}{2} & \frac{1}{2} & \frac{1}{2} & -\frac{1}{2} \\ 
\frac{1}{2} & -\frac{1}{2} & \frac{1}{2} & -\frac{1}{2} \\ 
\frac{1}{2} & -\frac{1}{2} & -\frac{1}{2} & \frac{1}{2} \\ 
\frac{1}{2} & -\frac{1}{2} & \frac{1}{2} & -\frac{1}{2}%
\end{array}%
\right) \text{ and }\rho \left( w\right) =\left( 
\begin{array}{cccc}
-\frac{1}{2} & \frac{1}{2} & \frac{1}{2} & -\frac{1}{2} \\ 
\frac{1}{2} & -\frac{1}{2} & -\frac{1}{2} & \frac{1}{2} \\ 
\frac{1}{2} & \frac{1}{2} & -\frac{1}{2} & -\frac{1}{2} \\ 
\frac{1}{2} & \frac{1}{2} & -\frac{1}{2} & -\frac{1}{2}%
\end{array}%
\right) 
\]%
are tripotent matrices.

For $z=\frac{1}{2}\left( f_{1}+f_{2}+\sqrt{2}f_{3}\right) $, the matrices 
\[
\varphi \left( z\right) =\left( 
\begin{array}{cccc}
0 & \frac{1}{2} & \frac{1}{2} & -\frac{\sqrt{2}}{2} \\ 
\frac{1}{2} & 0 & \frac{\sqrt{2}}{2} & -\frac{1}{2} \\ 
\frac{1}{2} & -\frac{\sqrt{2}}{2} & 0 & \frac{1}{2} \\ 
\frac{\sqrt{2}}{2} & -\frac{1}{2} & \frac{1}{2} & 0%
\end{array}%
\right) \text{ and }\rho \left( z\right) =\left( 
\begin{array}{cccc}
0 & \frac{1}{2} & \frac{1}{2} & -\frac{\sqrt{2}}{2} \\ 
\frac{1}{2} & 0 & -\frac{\sqrt{2}}{2} & \frac{1}{2} \\ 
\frac{1}{2} & \frac{\sqrt{2}}{2} & 0 & -\frac{1}{2} \\ 
\frac{\sqrt{2}}{2} & \frac{1}{2} & -\frac{1}{2} & 0%
\end{array}%
\right) 
\]%
are nilpotent matrices with $2$ as a nilpotency index.

ii) If we consider the split quaternion algebra $\mathbb{H}\left( 2,3\right) 
$ and the quaternion $q_{1}=\frac{1}{2}\left( 1+f_{1}+f_{2}+\frac{\sqrt{6}}{3%
}f_{3}\right) $, with $\mathbf{n}\left( q_{1}\right) =0$,we have that the
matrices 
\[
\varphi \left( q_{1}\right) =\left( 
\begin{array}{cccc}
\frac{1}{2} & 1 & \frac{3}{2} & -\sqrt{6} \\ 
\frac{1}{2} & \frac{1}{2} & \frac{\sqrt{6}}{2} & -\frac{3}{2} \\ 
\frac{1}{2} & -\frac{\sqrt{6}}{3} & \frac{1}{2} & 1 \\ 
\frac{\sqrt{6}}{6} & -\frac{1}{2} & \frac{1}{2} & \frac{1}{2}%
\end{array}%
\right) \text{ and }\rho \left( q_{1}\right) =\left( 
\begin{array}{cccc}
\frac{1}{2} & 1 & \frac{3}{2} & -\sqrt{6} \\ 
\frac{1}{2} & \frac{1}{2} & -\frac{\sqrt{6}}{2} & \frac{3}{2} \\ 
\frac{1}{2} & \frac{\sqrt{6}}{3} & \frac{1}{2} & -1 \\ 
\frac{\sqrt{6}}{6} & \frac{1}{2} & -\frac{1}{2} & \frac{1}{2}%
\end{array}%
\right) 
\]%
are idempotent. 

If we take $q_{2}=\frac{1}{2}\left( -1+f_{1}+f_{2}+\frac{\sqrt{6}}{3}%
f_{3}\right) $, with $\mathbf{n}\left( q_{2}\right) =0$,we have that the
matrices 
\[
\varphi \left( q_{1}\right) =\left( 
\begin{array}{cccc}
-\frac{1}{2} & 1 & \frac{3}{2} & -\sqrt{6} \\ 
\frac{1}{2} & -\frac{1}{2} & \frac{\sqrt{6}}{2} & -\frac{3}{2} \\ 
\frac{1}{2} & -\frac{\sqrt{6}}{3} & -\frac{1}{2} & 1 \\ 
\frac{\sqrt{6}}{6} & -\frac{1}{2} & \frac{1}{2} & -\frac{1}{2}%
\end{array}%
\right) \text{ and }\rho \left( q_{1}\right) =\left( 
\begin{array}{cccc}
-\frac{1}{2} & 1 & \frac{3}{2} & -\sqrt{6} \\ 
\frac{1}{2} & -\frac{1}{2} & -\frac{\sqrt{6}}{2} & \frac{3}{2} \\ 
\frac{1}{2} & \frac{\sqrt{6}}{3} & -\frac{1}{2} & -1 \\ 
\frac{\sqrt{6}}{6} & \frac{1}{2} & -\frac{1}{2} & -\frac{1}{2}%
\end{array}%
\right) 
\]%
are tripotent. 

If we take $q_{3}=\frac{1}{2}\left( 1+\sqrt{2}f_{1}+f_{2}+f_{3}\right) $,
the matrices

\[
\varphi \left( q_{3}\right) =\left( 
\begin{array}{cccc}
\frac{1}{2} & \sqrt{2} & \frac{3}{2} & -3 \\ 
\frac{\sqrt{2}}{2} & \frac{1}{2} & \frac{3}{2} & -\frac{3}{2} \\ 
\frac{1}{2} & -1 & \frac{1}{2} & \sqrt{2} \\ 
\frac{1}{2} & -\frac{1}{2} & \frac{\sqrt{2}}{2} & \frac{1}{2}%
\end{array}%
\right) \text{and }\rho \left( q_{3}\right) =\left( 
\begin{array}{cccc}
\frac{1}{2} & \sqrt{2} & \frac{3}{2} & -3 \\ 
\frac{\sqrt{2}}{2} & \frac{1}{2} & -\frac{3}{2} & \frac{3}{2} \\ 
\frac{1}{2} & 1 & \frac{1}{2} & -\sqrt{2} \\ 
\frac{1}{2} & \frac{1}{2} & -\frac{\sqrt{2}}{2} & \frac{1}{2}%
\end{array}%
\right) 
\]%
are idempotent. For $q_{4}=\frac{1}{2}\left( -1+\sqrt{2}f_{1}+f_{2}+f_{3}%
\right) $, we obtain the following tripotent matrices

\[
\varphi \left( q_{4}\right) =\left( 
\begin{array}{cccc}
-\frac{1}{2} & \sqrt{2} & \frac{3}{2} & -3 \\ 
\frac{\sqrt{2}}{2} & -\frac{1}{2} & \frac{3}{2} & -\frac{3}{2} \\ 
\frac{1}{2} & -1 & -\frac{1}{2} & \sqrt{2} \\ 
\frac{1}{2} & -\frac{1}{2} & \frac{\sqrt{2}}{2} & -\frac{1}{2}%
\end{array}%
\right) \text{and }\rho \left( q_{4}\right) =\left( 
\begin{array}{cccc}
-\frac{1}{2} & \sqrt{2} & \frac{3}{2} & -3 \\ 
\frac{\sqrt{2}}{2} & -\frac{1}{2} & -\frac{3}{2} & \frac{3}{2} \\ 
\frac{1}{2} & 1 & -\frac{1}{2} & -\sqrt{2} \\ 
\frac{1}{2} & \frac{1}{2} & -\frac{\sqrt{2}}{2} & -\frac{1}{2}%
\end{array}%
\right) \text{.}
\]

$\allowbreak $\textbf{Conclusions.} The  $k-$potent matrices have many
applications in manny fields of research as for example combinatorics and
graph theoty, control theory, etc. For this reason, we considered that a
procedure to obtain some example of $k-$potent matrices of order $4$ or $8$,
over the real field or over the field $\mathbb{Z}_{p}$, with $p$ a prime
number is very usefull. The connections with quaternions and octonins
allowed us to obtain such examples. For a further research, we will study
the possibility to obtain new procedures which allow us to obtain new
classes and examples of $k-$potent matrices.      

\[
\]

\textbf{References}%
\[
\]

[FB; 24] Flaut, C., Baias, A., \textit{Some Remarks Regarding Special
Elements in Algebras Obtained by the Cayley--Dickson Process over} $\mathbb{Z%
}_{p},$Axioms, 13(6)(2024), 351;
https://doi.org/10.3390/axioms13060351\smallskip

[Sc; 54] Schafer, R. D., \textit{On the algebras formed by the
Cayley-Dickson process,} Amer. J. Math., \textbf{76}(1954),
435-446.\smallskip

\textbf{[}Sc; 66\textbf{]} Schafer, R. D., \textit{An Introduction to
Nonassociative Algebras,} Academic Press, New-York, 1966.\smallskip

[Ti; 00] Tian, Y., \textit{Matrix Representatios of Octonions and their
Applications}, https://arxiv.org/pdf/math/0003166, 2000.%
\[
\]

Cristina FLAUT

Faculty of Mathematics and Computer Science, Ovidius University,

Bd. Mamaia 124, 900527, Constan\c{t}a, Rom\^{a}nia,

http://www.univ-ovidius.ro/math/

e-mail: cflaut@univ-ovidius.ro; cristina\_flaut@yahoo.com

\[
\]

Andreea BAIAS

PhD student at Doctoral School of Mathematics,

Ovidius University of Constan\c{t}a, Rom\^{a}nia,

e-mail: andreeatugui@yahoo.com

\end{document}